\documentclass[final,onefignum,onetabnum]{siamart190516}
\usepackage{amsmath}
\usepackage{latexsym}
\usepackage{amssymb}
\usepackage{algorithm,algorithmic}
\usepackage{color}
 \usepackage{mathptmx}  

\def\hatr{\begin{array}}\def\ear{\end{array}}
\def\bas{\begin{eqnarray*}}\def\eas{\end{eqnarray*}}
\def\bea{\begin{eqnarray}\label}\def\eea{\end{eqnarray}}
\def\beq{\begin{equation}\label}\def\eeq{\end{equation}}
\def\bde{\begin{description}}\def\ede{\end{description}}
\def\ben{\begin{enumerate}}\def\een{\end{enumerate}}
\def\bit{\begin{itemize}}\def\eit{\end{itemize}}
\def\ben{\begin{enumerate}}\def\een{\end{enumerate}}
\def\bc{\begin{center}}\def\ec{\end{center}}
\def\mat{\left[\begin{array}}\def\rix{\end{array}\right]}

\def\beq{\begin{equation}\label}\def\eeq{\end{equation}}
\def\bea{\begin{eqnarray}\label}\def\eea{\end{eqnarray}}

\newtheorem{remark}{Remark}

\def\R{\mathbb{R}}

\usepackage{algorithm}
\usepackage{algorithmic}
\usepackage{caption}
\usepackage{subcaption}
\usepackage{hyperref}

\title{A multidimensional principal component analysis via the c-product Golub-Kahan-SVD for classification and  face recognition}

\author{M. Hached \thanks{Univ. Lille, CNRS, UMR 8524 - Laboratoire Paul Painlevé, F-59000 Lille, France; {\tt email: mustapha.hached@univ-lille.fr}} \and K. Jbilou \thanks{Laboratoire LMPA, 50 rue F. Buisson, ULCO Calais, France; {\tt email: jbilou@univ-littoral.fr} } \and C. Koukouvinos \thanks{Department of Mathematics, National Technical University of Athens, Zografou, 15773, Athens, Greece; {\tt email: ckoukouv@math.ntua.gr}}\and M. Mitrouli \thanks{National and Kapodistrian University of Athens
Panepistimiopolis, GR-157 84, Athens, Greece; {\tt email: mmitroul@math.uoa.gr}} }

\begin {document}
\maketitle
\begin{abstract}

Face recognition and identification is a very important application in machine learning. Due to the increasing amount of available data,  traditional approaches based on matricization and matrix PCA methods can be difficult to implement. Moreover, the tensorial approaches are a natural choice, due to the mere structure of the databases, for example in the case of color images. Nevertheless,  even though various authors proposed factorization strategies for tensors, the size of the considered tensors can pose some serious issues.  When only a few features are needed to construct the projection space, there is no need to compute a SVD on the whole data. Two versions of the tensor Golub-Kahan algorithm are considered in this manuscript, as an alternative to the classical use of the tensor SVD which is based on truncated strategies. In this paper, we consider  the Tensor Tubal  Golub Kahan Principal Component Analysis method which purpose it to extract the main features of images using the tensor singular value decomposition (SVD) based on the tensor cosine product that uses the  discrete cosine transform.  This approach is applied for classification and face recognition and numerical tests show its effectiveness.\\

\textit{This article is dedicated to Constantin M. Petridi.}
 
\end{abstract}

\pagestyle{myheadings}
\thispagestyle{plain}
\markboth{M. Hached, K. Jbilou, C. Koukouvinos and M. Mitrouli}{A multidim. PCA via the c-product Golub-Kahan-SVD for classification and  face recognition.}


\section{Introduction} An important challenge in the last few years was the extraction of the main information in large datasets, measurements, observations that appear in signal and hyperspectral image processing, data mining, machine learning. Due to the increasing volume of data required by these applications, approximative low-rank matrix and tensor factorizations play a fundamental role in extracting latent components. The idea is to replace the initial large and maybe noisy and ill conditioned large scale original data by a lower dimensional approximate representation obtained via a matrix or multi-way array factorization or decomposition; see \cite{kilmer0,kolda1,lu,Martinez,elguide,beik2,ElIchi1,ElIchi2,brazell,vt2} for more details on recent work related to tensors and applications.  In the present work, we consider third order tensors that could be defined as three dimensional arrays of data. As our study is based on the cosine transform product, we limit this work to three-order tensors. 

The number of indices of a tensor is called modes or ways. Scalars can also be regarded as zero mode tensors, first mode tensors are vectors and matrices are second mode tensors. 

For a given 3-mode  tensor $ \mathcal {X}\in \mathbb{R}^{n_{1}\times n_{2}\times n_{3}}$, 
the notation $x_{{i_1},{i_2},i_3}$ stands for the element
 $\left(i_{1},i_{2},i_{3} \right) $ of the tensor $\mathcal {X}$. 

A fiber is defined by fixing all the indexes  except  one. A matrix column is a mode-1 fiber and a matrix row is a mode-2 fiber. Third-order tensors have column, row and tube fibers. An element $c\in \mathbb{R}^{1\times 1 \times n}$ is called a tubal-scalar or simply tube of length $n$. More details can be found in  \cite{kilmer0,kolda1}.

\section{Definitions and Notations}

		\subsection{Discrete Cosine Transformation}
	\noindent In this subsection we recall some definitions and properties of the discrete cosine transformation and the c-product of tensors. During recent years, many advances were made in order to establish a rigorous framework enabling the treatment of problems for which the data is stored in three-way tensors without having to resort to matricization \cite{lu,kilmer0}. One of the most important feature of such a framework is the definition of a tensor-tensor product as the t-product, based on the Fast Fourier Transform . For applications as image treatment, the tensor-tensor product based on the Discrete Cosine Transformation (DCT) has shown to be an interesting alternative to FFT. We now give some basic facts on the DCT and its associated tensor-tensor product. The DCT of a vector $v \in {\R}^n$ is defined by
	\begin{equation}
	\label{dft1}
	\tilde v= C_n v \in {\R}^n,
	\end{equation}
	where $C_n$ is the $n\times n$  discrete cosine
	transform matrix with entries 
	$$ (C_n)_{ij}=\sqrt{\frac{ 2-\delta_{i1}}{n}}\cos\left( \dfrac{(i-1)(2j-1)\pi}{2n}\right)\quad 1 \leq i,j \leq n $$
	with $\delta_{ij} $ is the Kronecker delta; see [ \cite{Jain} p. 150] for more details. It is  known that the matrix $C_n$ is orthogonal, \textit{ie} $C_n^TC_n=C_nC_n^T=I_n$; see \cite{Ng}. Furthermore, for  any  vector $v \in \mathbb{R}^n$, the   matrix vector multiplication $C_nv$     can be computed in $O(nlog(n))$ operations.  Moreover, Ng  and  al. \cite{Ng}  have shown that a certain class of Toeplitz-plus-Hankel matrices can be diagonalized by $C_n$. More precisely, we have 
	\begin{equation}\label{dft5}
	C_n \, {\tt th }(v)\, C_n^{-1}  = {\rm Diag}(\tilde v),\\
	\end{equation}
	where
	$$
	{\tt th}(v)=  \mathop{ {\underbrace{\left ( 
				\begin{array}{cccc}
				v_1 & v_2 & \ldots & v_n\\
				v_2 & v_1 & \ldots & v_3\\
				\vdots & \vdots  & \ldots & \vdots \\
				v_n & v_{n-1} & \ldots & v_1\\
				\end{array}
				\right )}}}\limits_{\tiny{\text{Toeplitz}}}  +  \mathop{ {\underbrace{\left( {\begin{array}{*{20}{c}} {{v_2}}&\ldots &  {{v_n}} &0\\
					\vdots &  \reflectbox{ $\ddots$}   & \reflectbox{ $\ddots$} &   {{v_{{n}}}}\\
					{{v_{n}}}&0& \ldots&  \vdots\\
					0&{{v_{{n}}}}&\ldots&  {{v_2}}
					\end{array}} \right)}}}\limits_{\tiny{\text{Hankel}}}   
	$$
	and $ {\rm Diag}(\tilde v)$ is the diagonal matrix whose $i$-th diagonal element is $ (\tilde v)_i$.

	\subsection{Definitions and properties of the cosine product}
\medskip
	\noindent In this subsection, we briefly review some concepts and notations, which play a central role for the elaboration of the tensor global iterative methods based on the c-product; see \cite{kilmer5,reichel} for more details on the c-product.\\
	\medskip 	Let $\mathcal {A} \in \mathbb{R}^{n_{1}\times n_{2}\times n_{3}} $ be a real valued third-order tensor, then the operations {\tt mat}  and  its inverse {\tt ten} are defined by
	\begin{align*}
	{\tt mat}(\mathcal {A})=   \mathop{ {\underbrace{\left ( 
				\begin{array}{cccc}
				A_1 & A_2 & \ldots & A_n\\
				A_2 & A_1 & \ldots & A_3\\
				\vdots & \vdots  & \ldots & \vdots \\
				A_n & A_{n-1} & \ldots & A_1\\
				\end{array}
				\right )}}}\limits_{\tiny{\text{Block Toeplitz}}}  &+  \mathop{ {\underbrace{\left( {\begin{array}{*{20}{c}} {{A_2}}&\ldots &  {{A_n}} &0\\
					\vdots &  \reflectbox{ $\ddots$}   & \reflectbox{ $\ddots$} &   {{A_{{n}}}}\\
					{{A_{n}}}&0& \ldots&  \vdots\\
					0&{{A_{{n}}}}&\ldots&  {{A_2}}
					\end{array}} \right)}}}\limits_{\tiny{\text{Block  Hankel}}}  \in {\R}^{ n_1n_3 \times n_2n_3}
	\end{align*}				
	and the inverse operation denoted by  {\tt ten} is simply defined by 
	$$ {\tt ten}({\tt mat}(\mathcal {A}) ) =  \mathcal {A}.$$
	\noindent Let us denote $\widetilde {\mathcal{A}}$ the tensor obtained by applying the DCT on all the tubes of the tensor $\mathcal {A}$. This operation and its inverse are implemented in the Matlab by the commands ${\tt dct}$  and ${\tt idct}$ as 
	$$\widetilde {\mathcal{A}}= {\tt dct}(\mathcal {A},[ \;],3), \; {\rm and }\;\; {\tt idct} (\widetilde {\mathcal{A}}, [ \;],3)= \mathcal {A},$$
	where ${\tt idct}$ denotes the Inverse Discrete Cosine Transform.\\
\begin{remark}
		Notice that   the  tensor $\widetilde {\mathcal{A}}$ can be computed   by  using  the  3-mode product defined in   \cite{kolda1} as  follows:
		$$ \widetilde{{\mathcal{A} }}= {\mathcal{A} }\times_3 M    $$ 
		where  M is the  ${ n_{3}\times n_{3}}$  invertible matrix  given by  $$M=W^{-1}C_{n_3}(I+Z)$$ where $C_{n_3}$ denote de $ n_3\times n_3 $ Discrete Cosine   Transform  DCT matrix, $W= {\rm diag}(C_{n_3}(:,1))$ is the diagonal matrix made of the first column of the DCT matrix, Z is $n_3\times n_3 $ circulant upshift matrix which   can be computed in MATLAB using  $W={\rm diag}({\rm ones}(n_{3}-1,1),1)$ and I the $ n_3\times n_3 $ identity matrix   ; see \cite{kilmer5} for more details.\\
\end{remark}  
	Let ${\bf A}$ be the matrix 
	\begin{equation}\label{dft9}
	{\bf A}= \left (
	\begin{array}{cccc}
	{A}^{(1)}& &&\\
	& {A}^{(2)}&&\\
	&&\ddots&\\
	&&&{A}^{(n_3)}\\
	\end{array}
	\right)\in \mathbb{R}^{n_3 n_1\times n_3n_2 }
	\end{equation}
	where  the matrices ${A}^{(i)}$'s are the frontal slices of the tensor ${\widetilde {\mathcal{A}}}$. The  block   matrix ${\tt mat}(\mathcal {A})$ can also be block diagonalized by using the DCT matrix as follows
	\begin{equation}\label{dft8}
	(C_{n_3} \otimes I_{n_1})\, {\tt mat}(\mathcal {A})\, 	(C_{n_3}^{T} \otimes I_{n_2})={\bf A}
	\end{equation}
	\begin{definition}
		The \textbf{c-product}  of two tensors
		$\mathcal {A} \in \mathbb{R}^{n_{1}\times n_{2}\times n_{3}} $ and $\mathcal {B} \in \mathbb{R}^{n_{2}\times m\times n_{3}} $ is the  ${n_{1}\times m\times n_{3}}$ tensor  defined by:	
		$$\mathcal {A} \star_c \mathcal {B}={\tt ten}({\tt mat}(\mathcal {A}){\tt mat}(\mathcal {B}) ).$$
	\end{definition}
	Notice that from the relation \ref{dft9}, we can show that the   product $\mathcal {C}= \mathcal {A} \star_c \mathcal {B}$ is equivalent to $  {\bf C}= {\bf A}\,{\bf B}$. 
	The following algorithm allows us to compute, in an efficient way, the c-product of the tensors $\mathcal {A}$ and $\mathcal {B}$, see \cite{kilmer5}.\\
	
\begin{algorithm}[!h]
\caption{Computing the  c-product }\label{algo1}
\textbf{Inputs}: $\mathcal {A} \in \mathbb{R}^{n_{1}\times n_{2}\times n_{3}} $ and $\mathcal {B} \in \mathbb{R}^{n_{2}\times m\times n_{3}} $\\
\textbf{Output}: $\mathcal {C}= \mathcal {A} \star_c \mathcal {B}  \in \mathbb{R}^{n_{1}\times m \times n_{3}} $
\begin{enumerate}
\item Compute $ \widetilde{{\mathcal{A} }}= {\tt dct}(\mathcal {A},[ \;],3)$ and $\mathcal {\widetilde B}={\tt dct}(\mathcal {B},[\; ],3)$.
\item Compute each frontal slices of $\mathcal {\widetilde C}$ by
			$$C^{(i)}= A^{(i)} B^{(i)}  $$
\item Compute $  {{\mathcal{C} }}= {\tt idct}(\mathcal {\widetilde C},[ \;],3)$  .	 \end{enumerate}
\end{algorithm}

	\noindent Next, give some definitions and remarks on the c-product and related topics.
	\begin{definition} 
		The identity tensor $\mathcal{I}_{n_{1}n_{1}n_{3}} $ is the tensor such  that  each frontal slice  of   $\widetilde{{\mathcal{I} }}_{n_{1}n_{1}n_{3}}$  is the identity matrix $I_{n_1n_1}$ .\\
		An $n_{1}\times n_{1} \times n_{3}$ tensor $\mathcal{A}$ is said to be invertible if there exists a tensor $\mathcal{B}$ of order  $n_{1}\times n_{1} \times n_{3}$  such that
		$$\mathcal{A}  \star_c \mathcal{B}=\mathcal{I}_{ n_{1}  n_{1}  n_{3}} \qquad \text{and}\qquad \mathcal{B}  \star_c \mathcal{A}=\mathcal{I}_{ n_{1}  n_{1}  n_{3}}.$$
		In that case, we denote $\mathcal{B}=\mathcal{A}^{-1}$. It is clear that 	$\mathcal{A}$ is invertible if and only if   ${\tt mat }(\mathcal{A})$ is invertible.\\

\noindent 

The  inner scalar product is defined by
		$$\langle \mathcal{A}, \mathcal{B} \rangle = \displaystyle \sum_{i_1=1}^{n_1} \sum_{i_2=1}^{n_2}  \sum_{i_3=1}^{n_3} a_{i_1 i_2 i_3}b_{i_1 i_2 i_3}$$
		and its  corresponding norm  is given by
		$ \Vert \mathcal{A} \Vert_F=\displaystyle \sqrt{\langle  \mathcal{A} ,  \mathcal{A}  \rangle}.$\\
		An $n_{1}\times n_{1} \times n_{3}$ tensor  $\mathcal{Q}$  is said to be orthogonal if
		$\mathcal{Q}^{T}   \star_c  \mathcal{Q}=\mathcal{Q} \star_c \mathcal{Q}^{T}=\mathcal{I}_{ n_{1}  n_{1}  n_{3}}.$\\

	\end{definition}
	\medskip

\begin{remark}
Another interesting way for  computing the scalar product and the associated norm is as follows:
		$\langle \mathcal{A}, \mathcal{B} \rangle = \displaystyle \frac{1}{n_3}  \langle  {\bf A}, {\bf B} \rangle$ and $  \Vert \mathcal{A} \Vert_F= \displaystyle \frac{1}{\sqrt{n_3}} \Vert {{\bf A}} \Vert_F,$
		where the block diagonal matrix ${\bf A}$  is defined by \eqref{dft9}.
	\end{remark}
	 \begin{definition}
	 	A tensor is called f-diagonal if its frontal slices are diagonal matrices. It is called upper triangular if all its frontal slices are upper triangular. 
	 \end{definition}
 
 \medskip

\noindent Next we recall the Tensor Singular Value Decomposition  of a tensor; more details can be found in \cite{kilmer1}.	 

\begin{theorem}\label{theosvd1}
	Let $\mathcal{A}$ be an $n_1 \times n_2 \times n_3$ real-valued tensor. Then $\mathcal{A}$ can be factored as follows
	\begin{equation}
	\label{svd1}
	\mathcal{A} = \mathcal{U} \star_c \mathcal{S} \star_c \mathcal{V}^T,
	\end{equation}
	where $\mathcal{U}$ and $\mathcal{V}$ are orthogonal tensors of order $(n_1,n_1,n_3)$ and $(n_2,n_2,n_3)$, respectively and $\mathcal{S}$ is an f-diagonal tensor of order $(n_1 \times n_2 \times n_3)$. This factorization is called Tensor Singular Value Decomposition (c-SVD) of the tensor $	\mathcal{A} $.
\end{theorem}

\medskip

 \begin{algorithm}[!h]
	\caption{The Tensor SVD (c-SVD) }\label{T_SVD}
	{\bf Input}: $\mathcal {A} \in \mathbb{R}^{n_{1}\times n_{2}\times n_{3}} $
	{\bf Output}: $\mathcal {U}$, $\mathcal {V}$ and  $\mathcal {S}$.
	\begin{enumerate}
		\item Compute $\mathcal {\widetilde A}={\tt dct}(\mathcal {A},[ ],3)$. 
		\item Compute each frontal slices of $\mathcal {\widetilde U}$, $\mathcal {\widetilde V}$ and $\mathcal {\widetilde S}$ from $\mathcal {\widetilde A}$ as follows
		\begin{enumerate}
		\item 	for $i=1,\ldots,n_3$
		 $$[\mathcal {\widetilde U}^{(i)},\mathcal {\widetilde S}^{(i)},\mathcal {\widetilde V}^{(i)}]=svd(\mathcal {\widetilde A}^{(i)})$$
    \item end for
    \end{enumerate}
		\item Compute $\mathcal {U}={\tt idct}(\mathcal {\widetilde U},[\,],3)$, 	$\mathcal {S}={\tt idct}(\mathcal {\widetilde S},[\,],3)$ and $\mathcal {V}={\tt idct}(\mathcal {\widetilde V},[\,],3)$.
	\end{enumerate}
\end{algorithm}

\begin{remark}
As for the t-product \cite{kilmer1},  we can show that if $	\mathcal{A} = \mathcal{U} \star_c \mathcal{S} \star_c \mathcal{V}^T$ is a c-SVD of the tensor $\mathcal{A}$, then we have
	\begin{equation}
	\label{svd2}
	\displaystyle \sum_{k=1}^{n_3} A_k = \left (	\displaystyle \sum_{k=1}^{n_3} U_k \right )\, \left (	\displaystyle \sum_{k=1}^{n_3} S_k \right )\,  \left (	\displaystyle \sum_{k=1}^{n_3} V_k^T \right ),
	\end{equation}
	where $A_k$, $U_k$, $S_k$ and $V_k$ are the frontal slices of the tensors $\mathcal{A}$, $\mathcal{U}$, $\mathcal{S}$ and $\mathcal{V}$, respectively, and
	\begin{equation}
	\label{svd3}
	\mathcal{A}= \displaystyle \sum_{i=1}^{\min (n_1,n_2)} \mathcal{U}(:,i,:)  \star_c \mathcal{S}(i,i,:)  \star_c \mathcal{V}(:,i,:)^T.
	\end{equation}	
\end{remark}

\medskip 

\begin{theorem}\label{tsvd}
	Let $	\mathcal{A} = 		 \mathcal{U} \star_c \mathcal{S} \star_c \mathcal{V}^T$ given by \eqref{svd1}, and define  for $ k \le min (n_1,n_2)$  the tensor
	\begin{equation}
	\label{svd4}
	\mathcal{A}_k=\displaystyle \sum_{i=1}^{k} \mathcal{U}(:,i,:)  \star_c \mathcal{S}(i,i,:)  \star_c \mathcal{V}(:,i,:)^T.
	\end{equation}
	Then 
	\begin{equation}
	\label{svd5}
	\mathcal{A}_k=arg \min_{\mathcal{X} \in \mathcal{M}} \Vert \mathcal{A}_k -\mathcal{A} \Vert_F,
	\end{equation}
	where $\mathcal{M}=\{ \mathcal{X} \star_c \mathcal{Y} ; \; \mathcal{X} \in {\R}^{n_1 \times k \times n_3}, \, \mathcal{Y} \in {\R}^{k \times n_2 \times n_3} \}$.
\end{theorem}
Note that when $n_3=1$ this theorem reduces to the well known  Eckart-Young theorem for matrices \cite{golub1}.
\medskip	
\begin{definition}{\bf The tensor tubal-rank}\\
Let $\mathcal{A}$ be an $n_1 \times n_2 \times n_3$ be a tensor and consider its  c-SVD  $\mathcal{A} = \mathcal{U} \star_c \mathcal{S} \star_c \mathcal{V}^T$. The  tensor tubal rank of $\mathcal{A}$, denoted as rank$_t$($\mathcal{A}$) is defined to be the number of non-zero tubes of the f-diagonal tensor $\mathcal{S}$, i.e., 
$$ {\rm rank}_t(\mathcal{A} )=  \# \{i, \mathcal{S}(i,i,:) \ne 0 \}.$$
	\end{definition}

\begin{definition}\label{mr}
The multi-rank of the tensor $\mathcal{A}$  is a vector $p \in \mathbb{R}^{n_3}$  with the $ i$-th element equal to the rank of the $i$-th frontal slice of  $\mathcal{\widetilde A}={\tt fft}(\mathcal{A},[],3)$, i.e.
$$p(i)=rank(A^{(i)}), \, i=1,\ldots,n_3.$$
\end{definition}

\medskip

\noindent The well known QR matrix decomposition can also be extended to the tensor case; see \cite{kilmer1}
\begin{theorem}
	Let $\mathcal{A}$ be a real-valued tensor of order $n_1 \times n_2 \times n_3$. Then $\mathcal{A}$ can be factored as follows
	\begin{equation}
	\label{qr1}
	\mathcal{A}= \mathcal{Q} \star_c \mathcal{R},
	\end{equation}
	where $\mathcal{Q}$ is an $n_1 \times n_1 \times n_3$ orthogonal tensor and $\mathcal{R}$ is an $n_1 \times n_1 \times n_3$ f-upper triangular tensor.
\end{theorem}

\section{Tensor Principal Component Analysis for face recognition}

Principle Component Analysis PCA is a widely used technique in image classification and face recognition.  Many approaches involve a conversion of color images to grayscale in order to reduce the training cost.  Nevertheless, for some applications, color an is important feature and tensor based approaches offer the possibility to take it into account. Moreover, especially in the case of facial recognition, it allows the treatment of enriched databases including for instance additional biometric information. But, one have to bear in mind that the computational cost is an important issue as the volume of data can be very large. We first recall some background facts on the matrix based approach.

\subsection{The matrix case}

One of the simplest and most effective PCA approaches used in face recognition systems is the so-called eigenface approach. This approach transforms faces into a small set of essential characteristics, eigenfaces, which are the main components of the initial set of learning images (training set). Recognition is done by projecting a test image in the eigenface subspace, after which the person is classified by comparing its position in eigenface space with the position of known individuals. The advantage of this approach over other face recognition strategies resides in its simplicity, speed and insensitivity to small or gradual changes on the face.

\medskip
\noindent The process is defined as follows:  Consider a set of training faces $I_1$, $I_2$, $\ldots$, $I_p$. All the face images have the same size: $n \times m$.
Each face $I_i$ is transformed into a vector $x_i$ using the operation $vec$:   $x_i=vec(I_i)$. These vectors are columns of the $nm\times p$ matrix $$X=[x_1,\ldots, x_p].$$  We compute the average image $\mu=\displaystyle  \frac{1}{p} \displaystyle \sum_{i=1}^p x_i.$  Set 
$\bar x_i=x_i - \mu$ and consider the new matrices $$ \bar X=[\bar x_1,\ldots,\bar x_p],\;\; {\rm  and} \; \; C=\bar X { \bar X}^T.$$  Notice that the $nm \times nm$ covariance matrix   $C=\bar X { \bar X}^T$ can be very large. Therefore, the computation of the  $nm$ eigenvalues and the corresponding eigenvectors (eigenfaces) can be very difficult. To circumvent this issue, we instead consider the smaller $p \times p$ matrix $L={ \bar X}^T \bar X $.

\medskip
\noindent
Let $v_i$ be an eigenvector of $L$ then $Lv_i={ \bar X}^T{ \bar X}v_i=\lambda_i v_i$ and $${\bar X} Lv_i= {\bar X} {\bar X}^T {\bar X}v_i= \lambda_i {\bar X} v_i,$$
which shows that ${\bar X} v_i$ is an eigenvector of the covariance matrix $C={\bar X} {\bar X}^T$.

\noindent
The $p$ eigenvectors of $L={\bar X}^T {\bar X}$ are then used to find the $p$ eigenvectors $u_i={\bar X}v_i$ of $C$ that form the eigenface space. 
We keep only $k$ eigenvectors corresponding to the largest $k$ eigenvalues (eigenfaces corresponding to small eigenvalues can be omitted, as they explain only a small part of characteristic features of the faces.)

\noindent The next step consists in projecting each image of the training sample onto the eigenface space spanned by the orthogonal vectors $u_1,\ldots,u_k$:
$$\mathcal{U}_k=span\{u_1,\ldots,u_k\},\; {\rm with} \; { U}_k=[u_1,\ldots,u_k]$$
The matrix $U_kU_k^T$ is an orthogonal projector onto the subspace $\mathcal{U}_k$. A face image can be projected onto this face space as
$y_i ={U}_k^T (x_i-\mu).$ 

We now give the steps of an image classification process based on this approach:

\noindent Let $x=vec(I)$ be a test vector-image and project it onto the face space to get 
$y= {U}_k^T(x-\mu).$ Notice that the reconstructed image is given by
$$x^r= {\widetilde U}_k y +\mu.$$
Compute the Euclidean distance $$\epsilon_i= \Vert y - y_i \Vert,\, i=1,\ldots,k.$$
A face is classified as belonging to the class $l$ when the minimum $l$ is below some chosen threshold $\theta$ 
Set 
$$ \theta = \frac{1}{2} \displaystyle \max_{i,j} \Vert y_i - y_j \Vert,\; i,j=1,\ldots,k,$$  and let $\epsilon$ be the distance between the original test image $x$ ans its reconstructed image $x^r$:
 $ \epsilon = \Vert x-x^r \Vert$. Then 
\begin{itemize}
\item if $\epsilon \ge \theta$, then the input image is not even a face image and not recognized.

\item If $\epsilon<  \theta$ and $\epsilon_i \ge \theta$ for all  $i$ then the input image is a face image but it is an unknown image face.

\item If $\epsilon < \theta$  and $\epsilon_i <  \theta$ for all  $i$ then the input images are the individual face images associated with the class vector $x_i$.
\end{itemize}

\medskip
\noindent
We now give some basic facts on the relation between the singular value decomposition (SVD) and PCA in this context:

\noindent 
Consider the Singular Value Decomposition of the matrix $A$ as
$${\bar X} =U \Sigma V^T=\displaystyle \sum_{i=1}^p\sigma_i u_i v_i^T$$
where $U$ and $V$ are orthonormal matrices of sizes $nm$ and $p$ respectively. 
The singular values $\sigma_i$ are the square roots of the eigenvalues of the matrix $L={\bar X} ^T {\bar X}$, the $u_i$'s are the left vectors and the $v_i's$ are the right vectors. 
We have $$L={\bar X} ^T {\bar X} =V \Delta V^T;\;\; \Delta=diag(\sigma_1^2,\ldots,\sigma_p^2)$$
which is is the eigendecomposition of the matrix $L$ and
$$C={\bar X}  {\bar X} ^T=UD U^T; \; D=diag(\sigma_1^2,\ldots,\sigma_p^2,0,\ldots,0).$$
In the  PCA  method, the projected eigenface space is then  generated by the first $u_1,\ldots,u_k$ columns of the unitary matrix $U$ derived from the SVD decomposition of the matrix ${\bar X} $.\\

\medskip
As  only a small number $k$ of the largest singular values are needed in PCA, we can use the well known Golub-Kahan algorithm to compute these wanted singular values and the corresponding singular vectors to define the projected subspace.

\noindent In the next section, we explain how the SVD based PCA can be extended to tensors and propose an algorithm for facial recognition in this context.

\section{The tensor Golub-Kahan method}

As explained in the previous section, it is important to take into account the potentially large size of data, especially for the training process. The Golub Kahan bidiagonalization algorithm can be extended to the tensor context, especially in its c-tubal form.

\subsection{The Tensor c-global Golub Kahan algorithm}
Let $\mathcal{A} \in \mathbb{R}^{n_{1}\times n_{2}\times n_{3}}$  be a tensor ans $s\geq 1$ an integer. 
The Tensor c-global Golub Kahan bidiagonalization algorithm (associated to the c-product)  is defined as follows\\

\begin{algorithm}[h!]
	\caption{The Tensor Global Golub-Kahan algorithm (TGGKA)}\label{GGKA}
	\begin{enumerate}
		\item Choose a  tensor ${\mathcal V}_1 \in \mathbb{R}^{n_{2} \times s \times n_{3}}$ such that$\Vert \mathcal {V}_1 \Vert_F=1$  and set $\beta_0=0$.
		\item For $i=1,2,\ldots,k$
		\begin{enumerate}
			\item $\mathcal {U}_i=\mathcal {A} \star_c  \mathcal {V}_i-\beta_{i-1}\mathcal {U}_{i-1}$,
			\item $\alpha_i=\Vert \mathcal {U}_i \Vert_F$,
			\item $\mathcal {U}_i=\mathcal {U}_i/\alpha_i$,
			\item $\mathcal {V}_{i+1}=\mathcal {A}^T \star_c \mathcal {U}_i-\alpha_i \mathcal {V}_i$,
			\item $\beta_i=\Vert \mathcal {V}_{i+1} \Vert_F$.
			\item $\mathcal {V}_{i+1}=\mathcal {V}_{i+1}/\beta_i$.
		\end{enumerate}
		End
	\end{enumerate}
\end{algorithm}

\noindent Let $C_k$ be the $k \times k$ upper bidiagonal matrix defined by 
\begin{equation}\label{GCmat}
C_k=\begin{bmatrix}
\alpha_1&\beta_1\\
&\alpha_2 & \beta_2\\
&&\ddots&\ddots\\
& & &\alpha_{k-1}&\beta_{k-1}\\
&&&&\alpha_k
\end{bmatrix}.
\end{equation}
\medskip 

\noindent  Let  $\mathbb{V}_{k}  $ and $\mathcal{A}\star_c\mathbb{V}_{k}  $ be the $(n_{2}\times (sk)\times p)$ and   $(n_{1}\times (sk)\times n_{3})$ tensors with frontal slices $\mathcal{V}_{1},\ldots,\mathcal{V}_{k}$ and  $\mathcal{A}\star_c\mathcal{V}_{1},\ldots,\mathcal{A}\star_c\mathcal{V}_{k}$, respectively, and let  $\mathbb{U}_{k}  $ and $\mathcal{A}^T\star_c\mathbb{U}_{k}  $ be the $(n_{1}\times (sk)\times n_{3})$ and $(n_{2}\times (sk)\times n_{3})$  tensors with frontal slices $\mathcal{U}_{1},\ldots,\mathcal{U}_{k}$ and  $\mathcal{A}^T\star_c\mathcal{U}_{1},\ldots,\mathcal{A}^T\star_c\mathcal{U}_{k}$, respectively. We set  
\begin{align}
\label{ev12}
\mathbb{V}_{k}:&=\left[  \mathcal{V}_{1},\ldots,\mathcal{V}_{k}\right], \;\;\; {\rm and}\;\;\; \mathcal{A}\star_c\mathbb{V}_{k}:=[\mathcal{A}\star_c\mathcal{V}_{1},\ldots,\mathcal{A}\star_c\mathcal{V}_{k}],\\
\mathbb{U}_{k}:&=\left[  \mathcal{U}_{1},\ldots,\mathcal{U}_{k}\right], \;\;\; {\rm and} \;\;\; \mathcal{A}^T\star_c\mathbb{U}_{k}:=[\mathcal{A}^T\star_c\mathcal{U}_{1},\ldots,\mathcal{A}^T\star_c\mathcal{U}_{k}],
\end{align}
with 
\[
{\widetilde C}_k^T=\begin{bmatrix}
C_k^T\\
\beta_{k}e_k^T
\end{bmatrix}\in\mathbb{R}^{(k+1)\times k},\; e_k^T =(0,0,\ldots,0,1)^T.
\]

\noindent  Then, we have the following results \cite{ElIchi1}\\

\begin{proposition}\label{proptggkb} 
	The tensors produced by the tensor c-global Golub-Kahan algorithm satisfy the following relations
	\begin{eqnarray} \label{equa20}
	\mathcal {A} \star_c \mathbb{V}_k& = & \mathbb{U}_k\circledast{ C}_k, \\	
	\mathcal{A}^{T} \star_c \mathbb{U}_{k}& = & \mathbb{V}_{k+1}  \circledast {\widetilde  C}_k^T\\
		&= & \mathbb{V}_k \circledast{C}_k^T  + {\beta}_{k}  \left[  \mathcal{O}_{n\times s\times p},\ldots,\mathcal{O}_{n_{1}\times s\times n_{3}},\mathcal{V}_{k+1}\right] ,  
	\end{eqnarray}
\end{proposition} 
where the product  $\circledast$ is defined by: $$\mathbb{U}_{k}\circledast y=\sum_{j=1}^{k} {y}_{j}\mathcal{  {V}}_{j},\; y= (y_1,\ldots,y_m)^T\in \mathbb{R}^k.$$
We set the following notation:
\begin{equation*}
\mathbb{U}_{k}\circledast {C}_k=\left[   \mathbb{U}_k\circledast C_{k}^1 ,\ldots,\mathcal{U}_{k}\circledast C_{k}^k \right],
\end{equation*}
where $ C_{k}^i$  is the $i$-th column of the matrix $ C_{k}$.\\
We note that since the matrix $C_k$ is bidiagonal, $T_k=C_k^T C_k$ is symmetric and tridiagonal and then Algorithm computes the same information as tensor global Lanczos  algorithm  applied to the symmetric  matrix $A^{\ast} \star_cA$.\\

 \subsection{Tensor tubal Golub-Kahan bidiagonalisation algorithm}
First, we introduce some new products that will be useful in this section.\\

\begin{definition} \label{scalarmatrix}\cite{ElIchi1}
	Let  ${\rm \bf a}  \in {\mathbb R}^{1\times 1 \times n_{3}} $ and $\mathcal{B}   \in {\mathbb R}^{n_{1}\times n_{2} \times n_{3}} $, the  tube fiber tensor product $({\rm \bf a}\divideontimes\mathcal{B})$ is an $(n_{1}\times n_{2} \times n_{3})$ tensor defined by 
	\begin{align*}
	{\rm \bf a}\divideontimes\mathcal{B} =\begin{pmatrix}
	{\rm \bf a}\star_c b(1,1,:)&\ldots&{\rm \bf a}\star_cb(1,n_{2},:) \\
	\vdots&\ddots&\vdots \\
	{\rm \bf a}\star_c b(n_{1},1,:)&\ldots&{\rm \bf a}\star_c b(n_{1},n_{2},:) \\
	\end{pmatrix}
	\end{align*}
\end{definition}
\begin{definition}\label{bloctens0}   
	Let  $\mathcal{A}   \in {\mathbb R}^{n_{1}\times m_{1} \times n_{3}} $, $\mathcal{B}\in {\mathbb R}^{n_{1}\times m_{2} \times n_{3}}$, $\mathcal{C}   \in {\mathbb R}^{n_{2}\times m_{1} \times n_{3}} $ and   $\mathcal{D}\in {\mathbb R}^{n_{2}\times m_{2} \times n_{3}}$ be  tensors. The block tensor
	$$\left[ {\begin{array}{*{20}{c}}
		{\mathcal{A}}&{\mathcal{B}} \\
		{\mathcal{C}}&{\mathcal{D}} \\
		\end{array}} \right]\in {\mathbb R}^{(n_{1}+n_2)\times (m_1+m_{2}) \times n_{3}} $$
	is defined by compositing the frontal slices of the four tensors.
\end{definition}

 \begin{definition}   Let	$\mathcal{A}=[\mathcal{A}_1,\ldots, \mathcal{A}_{n_2}]   \in {\mathbb R}^{n_{1}\times n_{2} \times n_{3}} $ 
 	where $\mathcal{A}_{i}\in {\mathbb R}^{n_{1}\times 1 \times n_{3}} $ ,     we denoted by {\tt TVect}( $\mathcal{A}$  )  the \textbf{tensor vectorization } operator : ${\mathbb R}^{n_{1}\times n_{2} \times n_{3}}\mapsto {\mathbb R}^{n_{1}n_2\times 1 \times n_{3}}$ 
 	obtained by superposing the laterals slices  $\mathcal{A}_i $ of $ \mathcal{A}$, for $i=1,\ldots,n_2$. In  others words, for  a tensor $\mathcal{A}=[\mathcal{A}_1,\ldots, \mathcal{A}_{n_2}]   \in {\mathbb R}^{n_{1}\times n_{2} \times n_{3}} $ where $\mathcal{A}_{i}\in {\mathbb R}^{n_{1}\times 1 \times n_{3}} $ , we have :  
		$${\tt TVect}(\mathcal{A})=\begin{pmatrix}
		\mathcal{A}_1\\
		\mathcal{A}_2 \\
		\vdots \\
        \mathcal{A}_{n_2}
		\end{pmatrix}\in {\mathbb R}^{n_{1}n_2\times 1 \times n_{3}} $$
\end{definition}
\begin{remark}
  The  {\tt TVect} operator transform a given  tensor on lateral slice.  	 Its  easy  to see that when we  take $p=1$, the  {\tt TVect} operator   coincides with the operation $vec$ which transform the matrix on vector.
\end{remark}
\begin{proposition}
	 Let $\mathcal{A}$ be  a tensor of  size   $ {\mathbb R}^{n_{1}\times n_{2} \times n_{3}} $ , we have
	  $$\Vert {{ \mathcal{A}}} \Vert_F=\Vert {\tt TVec}( \mathcal{A})\Vert_F$$
\end{proposition}
\begin{definition}
	Let	$\mathcal{A}=[\mathcal{A}_1,\ldots, \mathcal{A}_{n_2}]   \in {\mathbb R}^{n_{1}\times n_{2} \times n_{3}} $ 
	where $\mathcal{A}_{i}\in {\mathbb R}^{n_{1}\times 1 \times n_{3}} $. We define the range space of $\mathcal{A}$ denoted by ${\tt Range}(\mathcal{A}) $  as the c-linear span of the lateral
	slices of $\mathcal{A}$ 
	\begin{equation}\label{range}
	 {\tt Range}(\mathcal{A}) =\left\lbrace\mathcal{A}_1\star_c a(1,1,:)+\dots+\mathcal{A}_{n_2}\star_c a(n_2,n_2,:)| a(i,i,:)\in {\mathbb R}^{1\times 1 \times n_{3}}  \right\rbrace 
	\end{equation}
\end{definition}

 \begin{definition}\cite{ElIchi2}  	Let	$\mathcal{A}   \in {\mathbb R}^{n_{1}\times n_{2} \times n_{3}} $ and   $\mathcal{B}\in {\mathbb R}^{m_{1}\times m_{2} \times n_{3}}$, the \textbf{c-Kronecker product}   $\mathcal{A}\odot \mathcal{B}$  of $ \mathcal{A}$ and $\mathcal{B}$ is  the $n_{1}m_{1}\times n_{2}m_{2} \times n_{3} $ tensor in which the i-th frontal slice of his transformed tensor $\widetilde{(\mathcal{A}\odot\mathcal{B})}$ is given by	:
	\begin{align*}
	\widetilde{(\mathcal{A}\odot\mathcal{B})}_{i} &=     ( {A}^{(i)}\otimes   {B}^{(i)}), \;\; i=1,...,n_{3}\end{align*}
where $A^{(i)}$ and $B^{(i)}$ are the $i$-th frontal slices of the tensors $\widetilde {\mathcal A}={\tt dct} (\mathcal{A},[\,],3)$ and 	$\widetilde {\mathcal B}={\tt dct} (\mathcal{B},[],3)$, respectively.
\end{definition}

\medskip 
\noindent We introduce now  a normalization algorithm allowing us to decompose the  non-zero tensor  $\mathcal{C}\in {\mathbb R}^{n_{1}\times n_{2} \times n_{3}}$, such that: 
$$ \mathcal{C}={\rm \bf a}\divideontimes\mathcal{Q},\; {\rm with}\;\; \left\langle \mathcal{Q},\mathcal{Q}\right\rangle={\rm \bf e},$$ 
 where ${\rm \bf a}$ is an  invertible  tube fiber of size ${\rm \bf a}\in {\mathbb R}^{1\times 1 \times n_{3}}$ and $\mathcal{Q}\in {\mathbb R}^{n_{1}\times n_{2} \times n_{3}}$  and {\rm \bf e} is the  tube    fiber  ${\rm \bf e}\in {\mathbb R}^{1\times 1 \times n_{3}}$ defined by  ${\tt unfold}({\rm \bf e})  =(1,0,0\ldots,0)^T$. \\
\noindent This procedure is described in Algorithm \ref{normalization12}

\begin{algorithm}[h]
	\caption{Normalization algorithm (Normalize)}\label{normalization12}
	\begin{enumerate}
		\item 	{\bf Input.} $\mathcal{A}\in {\mathbb R}^{n_{1}\times n_{2} \times n_{3}} $ and a tolerance $ tol>0$.  
		\item 	{\bf Output.} The tensor $\mathcal{Q}$ and the tube fiber ${\rm \bf a}$.
		\item Set $\mathcal{\widetilde{Q}}=\text{{\tt dct}}(\mathcal{A},[],3) $   	\begin{enumerate}
			\item For $j=1,\ldots,n_3$
			
			\begin{enumerate}
				\item  $a_j				 	=||\widetilde{Q}^{(j)}||_F  $ 
				\item    if  $a_j>tol$, $\mathcal{\widetilde{Q}}^{(j)}=\displaystyle \frac{\mathcal{\widetilde{Q}}^{(j)} }{a_j} $	
			\item	else $\mathcal{\widetilde{Q}}_{j}=\text{rand}(n_{1},n_{2})$;   $a_j				 	=||\widetilde{Q}^{(j)}||_F  $  \\
			$\mathcal{\widetilde{Q}}^{(j)}=\displaystyle \frac{\mathcal{\widetilde{Q}}^{(j)} }{a_j}$	;  $a_j=0,$
			\end{enumerate}	
			\item End
		\end{enumerate}
		\item $ \mathcal{Q}  =\text{{\tt idct}}  (\mathcal{\widetilde{Q}} ,[],3)$, $ {\rm \bf a}  =\text{{\tt idct}}({\rm \bf a} ,[],3)$
		\item End
	\end{enumerate}
\end{algorithm}

\noindent Next, we give the Tensor Tube Global Golub-Kahan (TTGGKA) algorithm, see{ElIchi1}. Let $\mathcal{A} \in \mathbb{R}^{n_{1}\times n_{2}\times n_{3}}$  be a tensor and  let $s\geq 1$ be an integer. The Tensor Tube Global Golub-Kahan bidiagonalization process is defined as follows.
\begin{algorithm}[h!]
	\caption{The Tensor Tube Global Golub-Kahan algorithm (TTGGKA)}\label{TubeGGKA}
	\begin{enumerate}
		\item Choose a  tensor $\mathcal{V}_1  \in \mathbb{R}^{n_{2}\times  s \times n_{3}}$ such that $\langle  V_1, V_1 \rangle ={\rm \bf e}$  and set ${\rm \bf b}_0=0$.
		\item For $i=1,2,\ldots,k$
		\begin{enumerate}
			\item $\mathcal {U}_i=\mathcal {A} \star_c  \mathcal {V}_i-{\rm \bf b} _{i-1}\divideontimes\mathcal {U}_{i-1}$,
			\item $[\mathcal {U}_i, {\rm \bf a} _{i}]=Normalize(\mathcal {U}_{i})$.
			\item $\mathcal {V}_{i+1}=\mathcal {A}^T \star_c \mathcal {U}_i-{\rm \bf a} _{i}\divideontimes \mathcal {V}_i$,
			\item$[\mathcal {V}_{i+1}, {\rm \bf b} _{i}]=Normalize(\mathcal {V}_{i+1})$. 
		\end{enumerate}
		End
	\end{enumerate}
\end{algorithm}

\noindent Let $\mathcal{C}_k$ be the $k \times k\times n_{3}$ upper bidiagonal tensor (each frontal slice of $\mathcal{C}_k$  is a bidiagonal matrix) and $\mathcal{ \widetilde{C}}_k$ the  $k \times (k+1)\times n_{3}$  defined by 
\begin{equation}\label{TubeGCmat}
\mathcal C_k=\begin{bmatrix}
{\rm \bf a}_1&{\rm \bf b}_1\\
&{\rm \bf a}_2 & {\rm \bf b}_2\\
&&\ddots&\ddots\\
& & &{\rm \bf a}_{k-1}&{\rm \bf b}_{k-1}\\
&&&&{\rm \bf a}_k
\end{bmatrix},\; {\rm and}\; 
 \mathcal{ \widetilde{C}}_k=\begin{bmatrix}
{\rm \bf a}_1&{\rm \bf b}_1&\\
&{\rm \bf a}_2 & {\rm \bf b}_2&\\
&&\ddots&\ddots&\\
& & &{\rm \bf a}_{k-1}&{\rm \bf b}_{k-1}&\\
&&&&{\rm \bf a}_k&{\rm \bf b}_{k}
\end{bmatrix}.
\end{equation}
\medskip

\noindent  Let  $\mathbb{V}_{k}  $ and $\mathcal{A}\star_c\mathbb{V}_{k}  $ be the $(n_{2}\times (sk)\times n_{3})$ and   $(n_{1}\times (sk)\times n_{3})$ tensors with frontal slices $\mathcal{V}_{1},\ldots,\mathcal{V}_{k}$ and  $\mathcal{A}\star_c\mathcal{V}_{1},\ldots,\mathcal{A}\star_c\mathcal{V}_{k}$, respectively, and let  $\mathbb{U}_{k}  $ and $\mathcal{A}^T\star_c\mathbb{U}_{k}  $ be the $(n_{1}\times (sk)\times n_{3})$ and $(n_{2}\times (sk)\times n_{3})$  tensors with frontal slices $\mathcal{U}_{1},\ldots,\mathcal{U}_{k}$ and  $\mathcal{A}^T\star_c\mathcal{U}_{1},\ldots,\mathcal{A}^T\star_c\mathcal{U}_{k}$, respectively. We set  
\begin{align}
\label{tubeev12}
\mathbb{V}_{k}:&=\left[  \mathcal{V}_{1},\ldots,\mathcal{V}_{k}\right], \;\;\; {\rm and}\;\;\; \mathcal{A}\star_c\mathbb{V}_{k}:=[\mathcal{A}\star_c\mathcal{V}_{1},\ldots,\mathcal{A}\star_c\mathcal{V}_{k}],\\
\mathbb{U}_{k}:&=\left[  \mathcal{U}_{1},\ldots,\mathcal{U}_{k}\right], \;\;\; {\rm and} \;\;\; \mathcal{A}^T\star_c\mathbb{U}_{k}:=[\mathcal{A}^T\star_c\mathcal{U}_{1},\ldots,\mathcal{A}^T\star_c\mathcal{U}_{k}],
\end{align}

\noindent  Then, we have the following results\\

\begin{proposition}\label{tubeproptggkb} 
	The tensors produced by the tensor TTGGKA algorithm satisfy the following relations
	\begin{eqnarray} \label{tubeequa20}
	\mathcal {A} \star_c \mathbb{V}_k& = & \mathbb{U}_k\star_c(\mathcal{ C}_k\odot\mathcal{ I}_{ssn_{3}}), \\	
	\mathcal{A}^{T} \star_c \mathbb{U}_{k}& = & \mathbb{V}_{k+1}\star_c  ( \widetilde { \mathcal{ C}}_k^T\odot\mathcal{ I}_{ssn_{3}})\\
  	&= & \mathbb{V}_k \star_c  ( {   \mathcal{ C}}_k^T\odot\mathcal{ I}_{ssp})  +\mathcal{ V}_{k+1} \star_c(({\rm \bf b}_{k} \star_c {\rm \bf e}_{1,k,:})\odot   \mathcal{I}_{ssn_{3}})   ,  
	\end{eqnarray}
		where ${\rm \bf e}_{1,k,:} \in \mathbb{R}^{1\times k\times n_{3}}$ with $1$ in  the $(1,k,1)$ position and zeros in the other positions,   $\mathcal{I}_{ssn_3}\in \mathbb{R}^{s\times s\times n_{3}}$ the identity tensor and  ${\rm \bf b}_k$ is the fiber  tube in the $(k,k+1,:)$ position of the  tensor   $\mathcal{ \widetilde{C}}_k$.
\end{proposition}

\section{The tensor tubal PCA method}

\medskip
\noindent
In this section, we describe a  tensor-SVD based PCA method for order 3 tensors which naturally arise in problems involving images such as facial recognition. As for the matrix case, we consider a set of $N$ training images, each of one being encoded as $n_{1}\times n_{2} \times n_{3}$ real tensors $\mathcal{I}_{i}$, $1 \leq i \leq N$. In the case of RGB images, each frontal slice would contain the encoding for each color layer ($n_{3}=3$) but in order to be able to store additional features, the case $n_{3}>3$ could be contemplated.

Let us consider one training image $\mathcal{I}_{i_0}$. Each one of the $n_{3}$ frontal slices $I_{i_0}^{(j)}$ of $\mathcal{I}_{i_0}$ is resized into a column vector $vec(I_{i_0}^{(j)})$ of length $L=n_{1} \times n_{2}$ and we form a $L\times 1 \times n_{3} $ tensor $\mathcal{X}_{i_0}$ defined by $\mathcal{X}_{i_0}(:,:,j)=vec(I_{i_0}^{(j)})$. Applying this procedure to each training image, we obtain $N$ tensors $\mathcal{X}_i$ of size $L\times 1 \times n_{3}$. The average image tensor is defined as $\displaystyle{\bar{\mathcal{X}}=\frac{1}{N}\sum_{i=1}^{N}\mathcal{X}_i}$ and we define the $L\times N \times n_{3}$ training tensor $\displaystyle{\mathcal{X}=[\bar{\mathcal{X}_1},\dots,\bar{\mathcal{X}_N}]}$, where $\bar{\mathcal{X}_i}=\mathcal{X}_i-\bar{\mathcal{X}}$.

\medskip
\noindent

Let us now consider the c-SVD decomposition $\mathcal{X}=\mathcal{U}*_c\mathcal{S}*_c \mathcal{V}^T$ of $\mathcal{X}$, where $\mathcal{U}$ and $\mathcal{V}$ are orthogonal tensors of size $L\times L \times n_{3}$ and $N\times N \times n_{3}$ respectively and $S$ is a f-diagonal tensor of size $L\times N \times n_{3}$. 

\medskip
\noindent
In the matrix context, it is known that just a few singular values suffice  to capture the main features of an image, therefore, applying this idea to each one of the three color layers, an RGB image can be approximated by a low tubal rank tensor. Let us consider an image tensor $\mathcal{S}\in \mathbb{R}^{n_1\times n_2 \times n_3}$ and its c-SVD decomposition $\mathcal{S}=\mathcal{U} \star_c \mathcal{S} \star_c \mathcal{V}^T$. Choosing an integer $r$ such as $r\leq \mbox{min}(n_1,n_2)$, we can approximate $\mathcal{S}$ by the $r$ tubal rank tensor   $$\displaystyle{\mathcal{S}_r \approx\sum_{i=1}^{r}\mathcal{U}(:,i,:)*_c\mathcal{S}(i,i,:)*_c \mathcal{V}(:,i,:)^T}.$$

\begin{figure}[htbp]
  \centering
  \includegraphics[width=14cm,height=9cm]{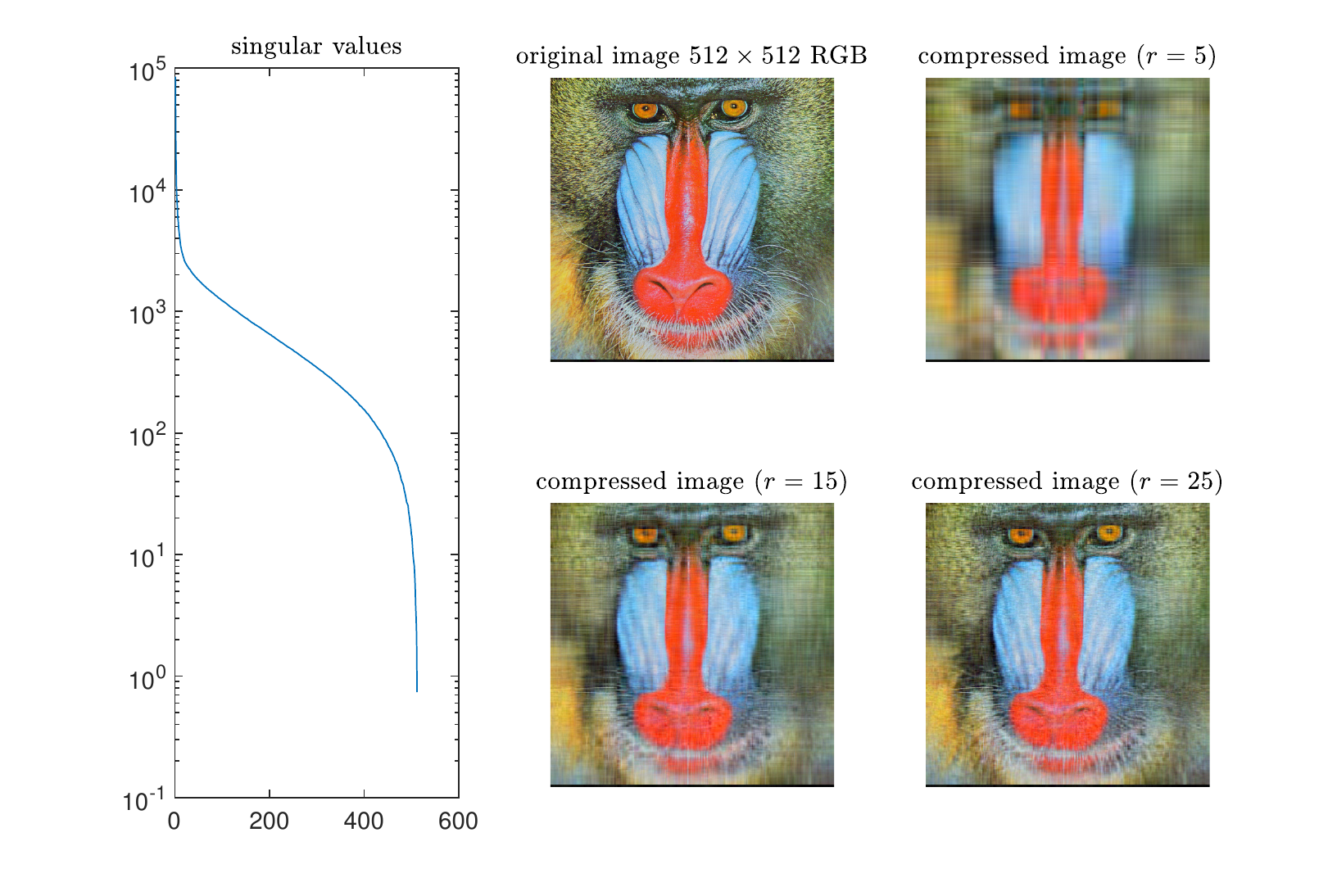}
  \caption{Image compression}
  \label{Fig1}
  
\end{figure}

In Figure \ref{Fig1}, we represented a $512 \times 512$ RGB image and the images obtained for various truncation indices. On the left part, we plotted the singular values of one color layer of the RGB tensor (the exact same behaviour is observed on the two other layers). The rapid decrease of the singular values explain the good quality of compressed images even for small truncation indices. 

\medskip

\noindent
Applying this idea to our problem, we want to be able to obtain truncated tensor SVD's of the training tensor $\displaystyle{\mathcal{X}}$, without needing to compute the whole c-SVD. After $k$ iterations of the TTGGKA algorithm (for the case $s=1$), we obtain three tensors $\mathbb{U}_k \in \mathbb{R}^{n_1 \times k \times n_3}$, $\mathbb{V}_{k+1} \in \mathbb{R}^{n_2 \times (k+1) \times n_3}$ and $\widetilde{\mathcal{C}}_k \in \mathbb{R}^{(k \times (k+1) \times n_3}$ as defined in \ref{tubeequa20} such as $$\mathcal{A}^{T} \star_c \mathbb{U}_{k} =  \mathbb{V}_{k+1}\star_c  \widetilde{\mathcal{C}}_k^T.$$

\begin{figure}
    \centering
    \includegraphics[width=9cm]{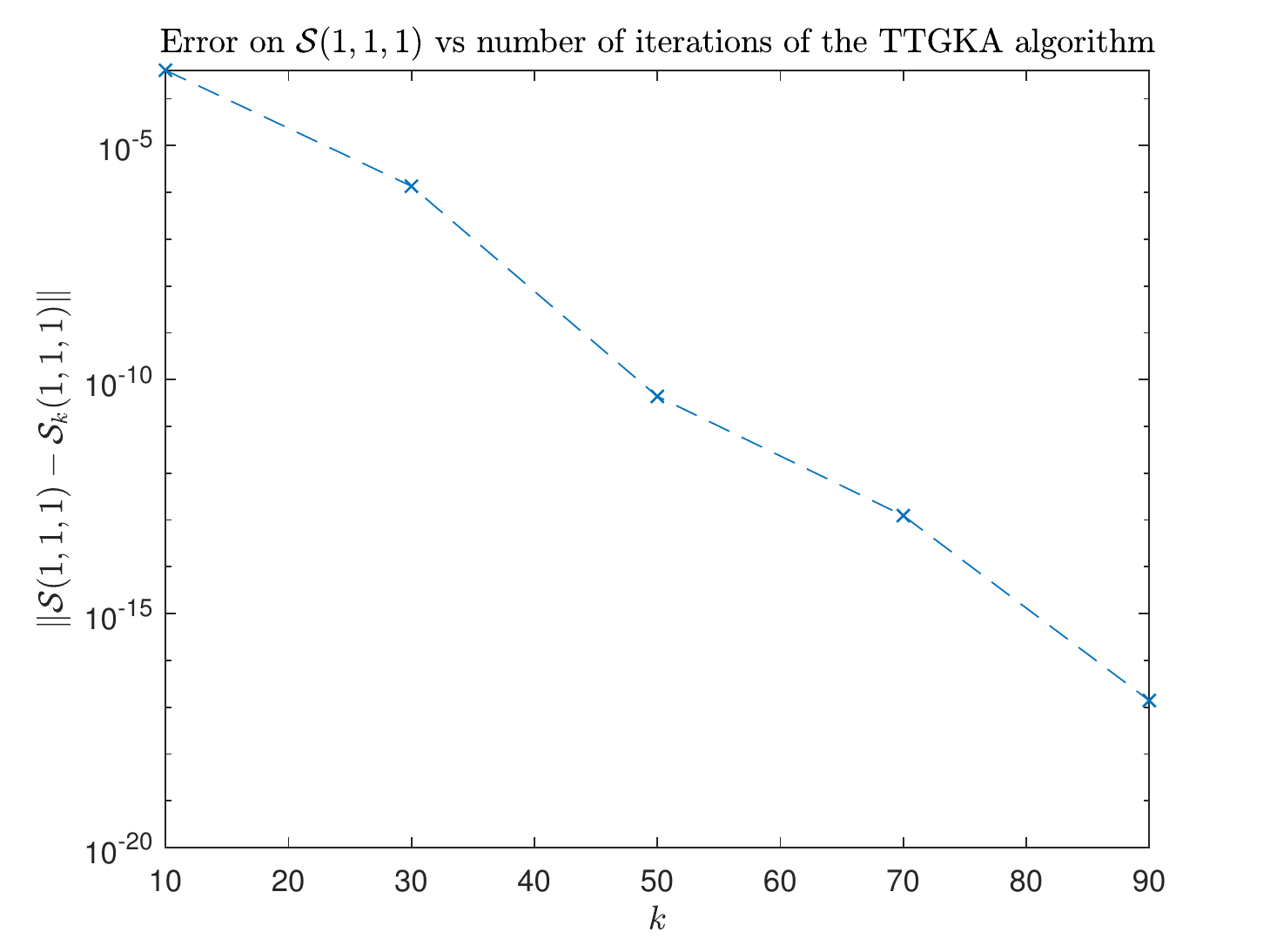}
    \caption{$\|\Sigma(1,1,1)-\mathcal{S}(1,1,1) \|$ \textit{vs} number of TTGGKA iteration $k$}
    \label{Fig2}
        \end{figure}

\noindent
Let $\displaystyle{\widetilde{\mathcal{C}}_k=\Phi \star_c \Sigma \star_c \Psi}$ the c-SVD of $\displaystyle{\widetilde{\mathcal{C}}_k}$, noticing that $\displaystyle{\widetilde{\mathcal{C}}_k \in \mathbb{R}^{k \times (k+1) \times n_{3}}}$ is much smaller than $\displaystyle{\bar{\mathcal{X}}}$. Then first tubal singular values and the left tubal singular tensors of $\displaystyle{\bar{\mathcal{X}}}$ are given by  $\displaystyle{\Sigma(i,i,:)}$ and $\displaystyle{\mathcal{U}_k\star_c \Phi(:,i,:)}$ respectively, for $i\leq k$, see \cite{kilmer0} for more details. 
\medskip
In order to illustrate the ability to approximate the first singular elements of a tensor using the TTGGKA algorithm, we considered a $900 \times 900 \times 3$ real tensor $\mathcal{A}$ which frontal slices were matrices generated by a finite difference discretization method of differential operators. On Figure \ref{Fig2}, we displayed the error on the first diagonal coefficient of the first frontal ${\mathcal S}(1,1,1)$ in function of the number of iteration of the Tensor Tube Golub Kahan algorithm, where $\displaystyle{\mathcal{A}=\mathcal{U}\star_c \mathcal{S} \star_c \mathcal{V}^T}$ is the c-SVD of $\mathcal{A}$.

\noindent     
In Table \ref{Tab1}, we reported on the errors on the tensor Frobenius norms of the singular tubes in function of the number $k$ of  the Tensor Tube Golub Kahan algorithm.
\medskip
\begin{table}[h!]
\begin{center}
\begin{tabular}{ l | c c c c }
  & $k=10$ & $k=30$ & $k=50$ &$k=70$\\
\hline \hline
$\mathcal{S}(1,1,:)$ & $3.6\times10^{-4}$ &$1.3\times10^{-5}$ & $5.1\times10^{-11}$ & $4.8\times10^{-17}$ \\ $\mathcal{S}(2,2,:)$ & $2.0\times10^{-3}$&$1.6\times10^{-6}$ & $5.2\times10^{-7}$ & $3.1\times10^{-8}$ \\
$\mathcal{S}(3,3,:)$ & $4.9\times10^{-3}$&$5.9\times10^{-4}$ & $2.3\times10^{-4}$ & $5.6\times10^{-8}$ \\
$\mathcal{S}(4,4,:)$ & $8.4\times10^{-3}$&$8.8\times10^{-4}$ & $1.5\times10^{-4}$ & $1.0\times10^{-8}$ \\
$\mathcal{S}(5,5,:)$ & $1.4\times10^{-2}$&$1.3\times10^{-3}$ & $2.7\times10^{-4}$ & $1.1\times10^{-8}$ \\
\hline
\end{tabular}
\end{center} 
\caption{$\|\mathcal{S}(i,i,:)-\Sigma(i,i,:)\|_F$ \textit{vs} $k$} 
\label{Tab1}
\end{table}  
\medskip

\noindent The same behaviour was observed on all the other frontal slices. This example illustrate the ability of the TTGKA algorithm for approximating the largest singular tubes.  
The projection space is generated by the lateral slices of the tensor  $\mathcal{P}=\mathbb{U}_k \star_c \Phi(:,1:k,:) \in \mathbb{R}^{n_1 \times i \times n_3}$ derived from the TTGGKA algorithm and the c-SVD decomposition of the bidiagonal tensor $\displaystyle{\widetilde{\mathcal{C}}_k}$ \textit{ie}, the c-linear span of first k lateral slices of $\mathcal{P}$, see\cite{kilmer1,kilmer0} for more details.  \\

\noindent The steps of the Tensor Tubal PCA algorithm for face recognition which finds the closest image in the training database for a given image $\mathcal{I}_0$ are summarized in Algorithm \ref{Tubepca}:\newpage

    \begin{algorithm}[h!]
    	\caption{The Tensor Tubal PCA algorithm (TTPCA )}\label{Tubepca}
\begin{enumerate}
\item \textbf{Inputs}  Training Image tensor $\mathcal{X}$ ($N$ images), mean image tensor $\bar{\mathcal{X}}$,Test image $\mathcal{I}_0$, index of truncation $r$, $k$=number of iterations of the TTGGKA algorithm ($k \geq r$). 
\item \textbf{Output} Closest image in the Training database.
\item Run $k$ iterations of the TTGGKA algorithm to obtain tensors $\mathbb{U}_k$ and $\widehat{\mathcal{C}_k}$
\item Compute $\left[ \Phi, \, \Sigma, \, \Psi \right]=$c-SVD$(\widetilde{\mathcal{C}_k})$
\item Compute the projection tensor $\mathcal{P}_r=\left[\mathcal{P}_r(:,1,:),\,\dots,\,\mathcal{P}_r(:,r,:)\right]$,\\ where $\mathcal{P}_r(:,i,:)=\mathbb{U}_k \star_c \Phi(:,i,:) \in \mathbb{R}^{n_1 \times 1 \times n_3}$
\item Compute the projected Training tensor $\hat{\mathcal{X}}_r=\mathcal{P}_r^T \star_c \mathcal{X}$ and projected centred test image $\hat{\mathcal{I}}_r=\mathcal{P}_r^T \star_c (\mathcal{I}-\bar{\mathcal{X}})$		
\item Find $i=\arg \min_{i=1,..,N} \|\hat{\mathcal{I}}_r\,-\,\hat{\mathcal{X}}_r(:,i,:)|_F$

    	\end{enumerate}
    \end{algorithm}

    
\medskip

\noindent In the next section, we consider image identification problems on various databases.

\section{Numerical tests}
In this section, we consider three examples of image identification. In the case of grayscale images, the global version of Golub Kahan was used to compute the dominant singular values in order to perform a PCA on the data. For the two other situations,  we used the Tensor Tubal PCA (TTPCA) method based on the Tube Global Golub Kahan (TTGGKA) algorithm in order to perform facial recognition on RGB images. The tests were performed with Matlab 2019a, on an Intel i5 laptop with 16 Go of memory.  We considered various truncation indices $r$ for which the recognition rates were computed. We also reported the CPU time for the training process (\textit{ie} the process of computating the eigenfaces). 

\subsection{Example 1}

In this example, we considered the MNIST database of handwritten digits \cite{MNIST}. The database contains two subsets of $28 \times 28$ grayscale images ( 60000 training images and 10000 test images). Each image was vectorized as a vector of length $28 \times 28=784$ and, following the process described in Section 3.1, we formed the training and the test matrices of sizes $784 \times 60000$ and $784 \times 10000$ respectively.

\begin{figure}[h!]
   \centering
   \includegraphics[width=8cm]{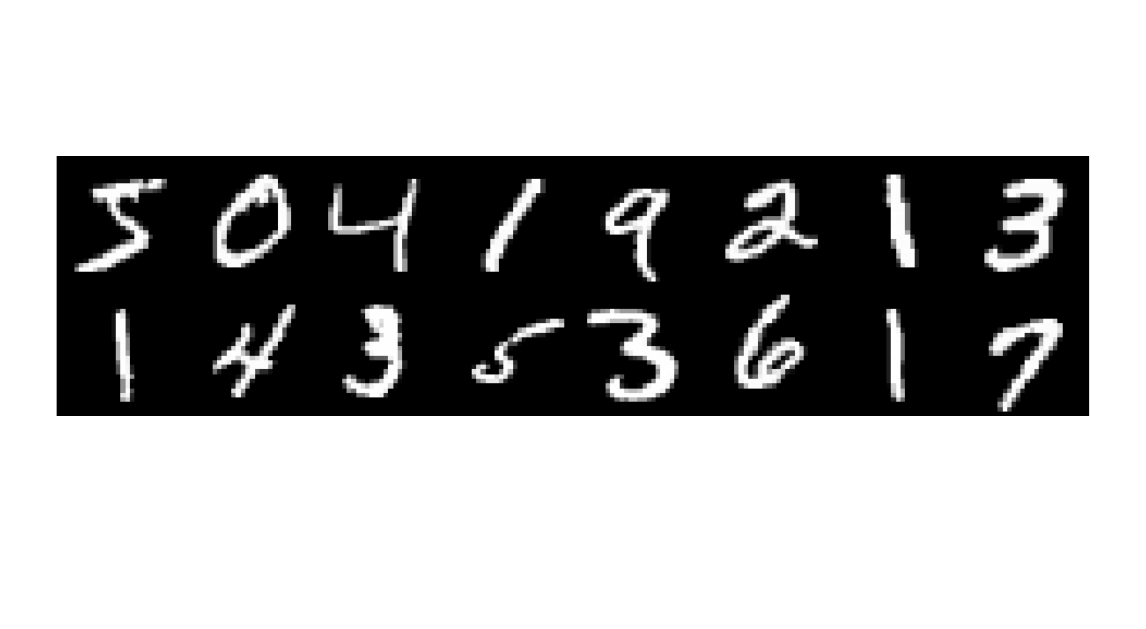}
   \caption{First 16 images of MNIST training subset}
   \label{Fig3}
   \end{figure}

\medskip
\noindent
Both matrices were centred by substracting the mean training image and the Golub Kahan algorithm was used to generate an approximation of $r$ dominant singular values $s_i$ and left singular vectors $u_i$, $i=1,\dots,r$.

\medskip

\noindent 
Let us denote $\mathcal{U}_r$ the subspace spanned by the columns of $U_r=\left[u_1,\dots,u_r\right]$. Let $t$ be a test image and $\hat{t}_r=U_r^Tt$ its projection onto $\mathcal{U}_r$. The closest image in the training dataset is determined by computing $$i=\arg \min_{i=1,..,60000} \| \hat{t}_r\,-\,\hat{X}_r(:,i)\|,$$
\noindent 
where $\hat{X}_r=\mathcal{U}_r^TX$.

\medskip

\begin{figure}[h!]
   \centering
   \includegraphics[width=8cm]{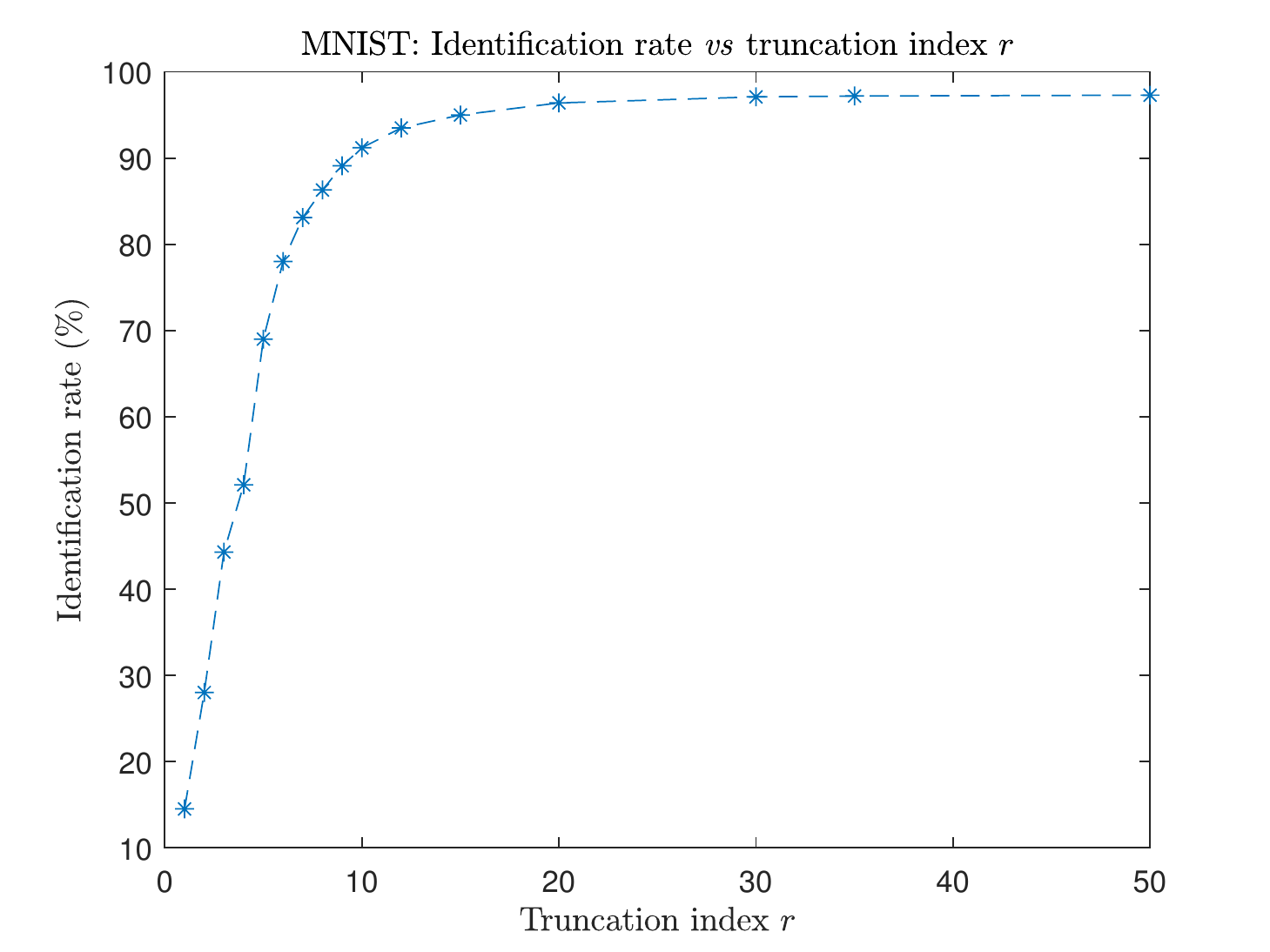}
   \caption{Identification rates for different truncation indices $r$}
   \label{Fig4}
   \end{figure}

\noindent
For various truncation indices $r$, we tested each image of the test subset and computed the recognition rate (\textit{ie} a test is successful if the digit is correctly identified). The results are plotted on Figure \ref{Fig4} and show that a good level of accuracy is obtained with only a few approximate singular values. Due to the large size of the training matrix, it validates the interest of computing only a few singular values with the Golub-Kahan algorithm.

\medskip

\subsection{Example 2}

In this example, we used the Georgia Tech database \href{http://www.anefian.com/research/face_reco.htm}{GTDB\_ crop} \cite{GTDB},  which contains 750 face images of 50 persons in different illumination conditions, facial expression and face orientation. The RGB JPEG images were resized to 100x100x3 tensors.

\begin{figure}[h!]
   \centering
\begin{tabular}{ccccc}
\includegraphics[width=1cm,height=1.35cm]{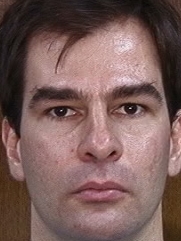}&
\includegraphics[width=1cm,height=1.35cm]{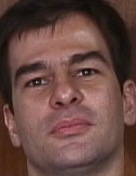}&
\includegraphics[width=1cm,height=1.35cm]{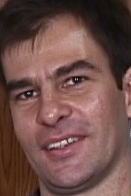}&
\includegraphics[width=1cm,height=1.35cm]{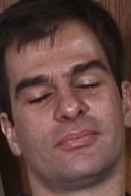}&
\includegraphics[width=1cm,height=1.35cm]{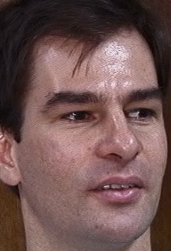}\\
\includegraphics[width=1cm,height=1.35cm]{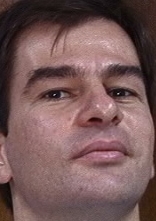}&
\includegraphics[width=1cm,height=1.35cm]{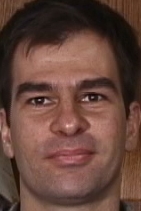}&
\includegraphics[width=1cm,height=1.35cm]{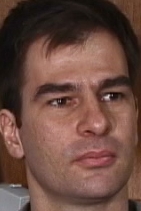}&
\includegraphics[width=1cm,height=1.35cm]{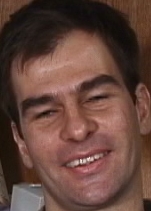}&
\includegraphics[width=1cm,height=1.35cm]{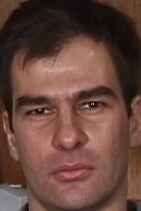}\\
\includegraphics[width=1cm,height=1.35cm]{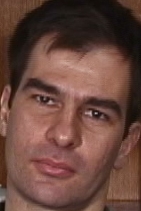}&
\includegraphics[width=1cm,height=1.35cm]{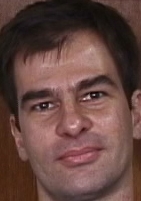}&
\includegraphics[width=1cm,height=1.35cm]{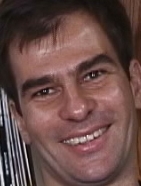}&
\includegraphics[width=1cm,height=1.35cm]{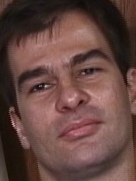}&
\includegraphics[width=1cm,height=1.35cm]{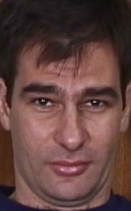}
\end{tabular}
    \caption{15 pictures of one individual in the database}
    \label{Fig5} 
\end{figure}

\noindent Each image file is coded as a $100 \times 100 \times 3$ tensor and transformed into a $10000 \times 1 \times 3$ tensor as explained in the previous section. We built the training and test tensors as follows: from 15 pictures of each person in the database, 5 pictures were randomly chosen and stored in the test folder and the 10 remaining pictures were used for the train tensor. Hence, the database was partitioned into two subsets containing 250 and 500 items respectively, at each iteration of the simulation.

\medskip

\begin{figure}[h!]
   \centering
   \includegraphics[width=8cm]{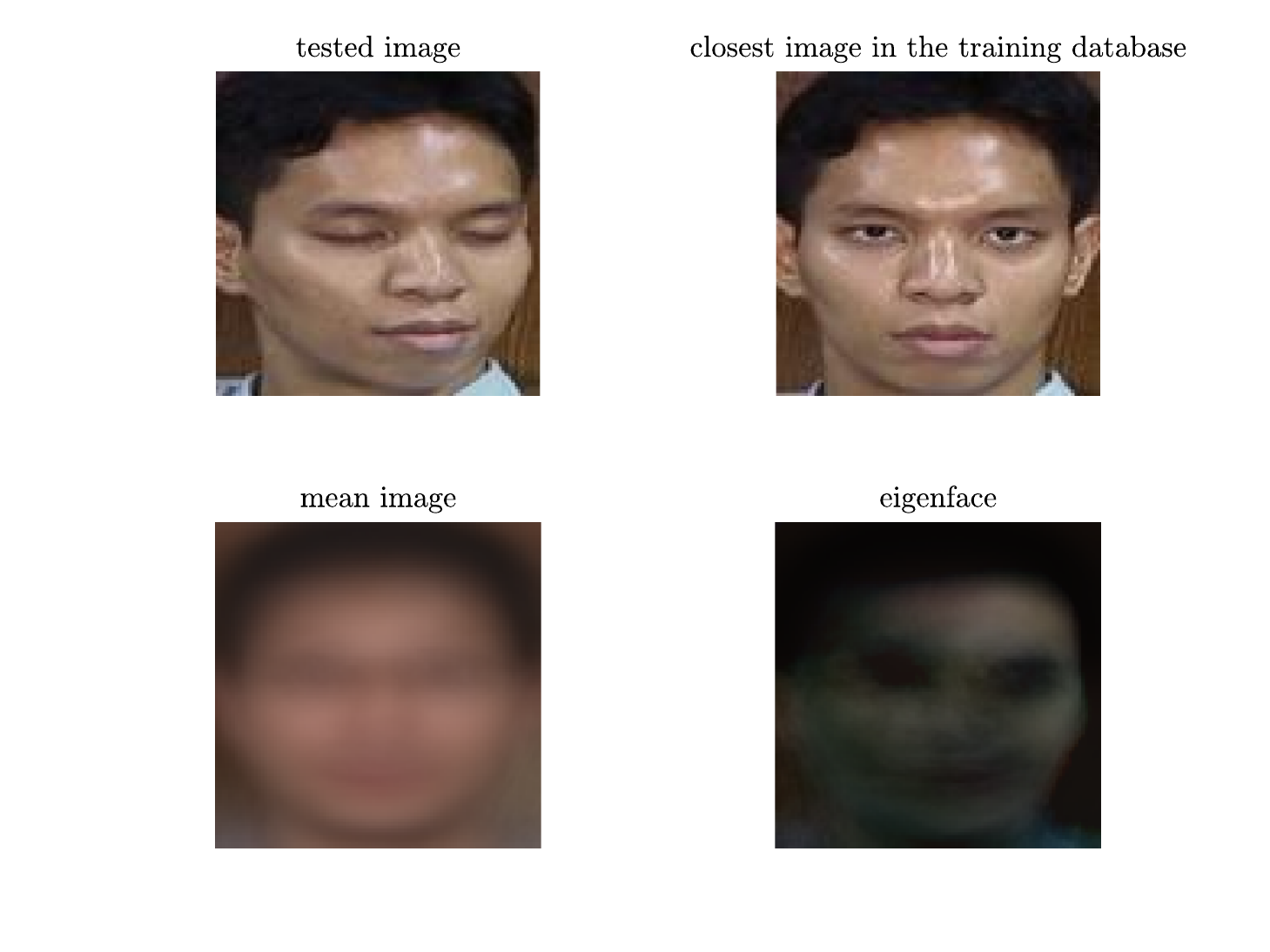}
   \caption{Test image, closest image, mean image and eigenface}
   \label{Fig6}
   \end{figure}

\noindent 
We applied the TTGGKA based algoritm \ref{Tubepca} for various truncation indices. In Figure \ref{Fig6}, we represented a test image (top left position), the closest image in the database (top right), the mean image of the training database (bottom left) and the eigenface associated to the test image (bottom right).

\begin{figure}[h!]
   \centering
   \includegraphics[width=8.5cm]{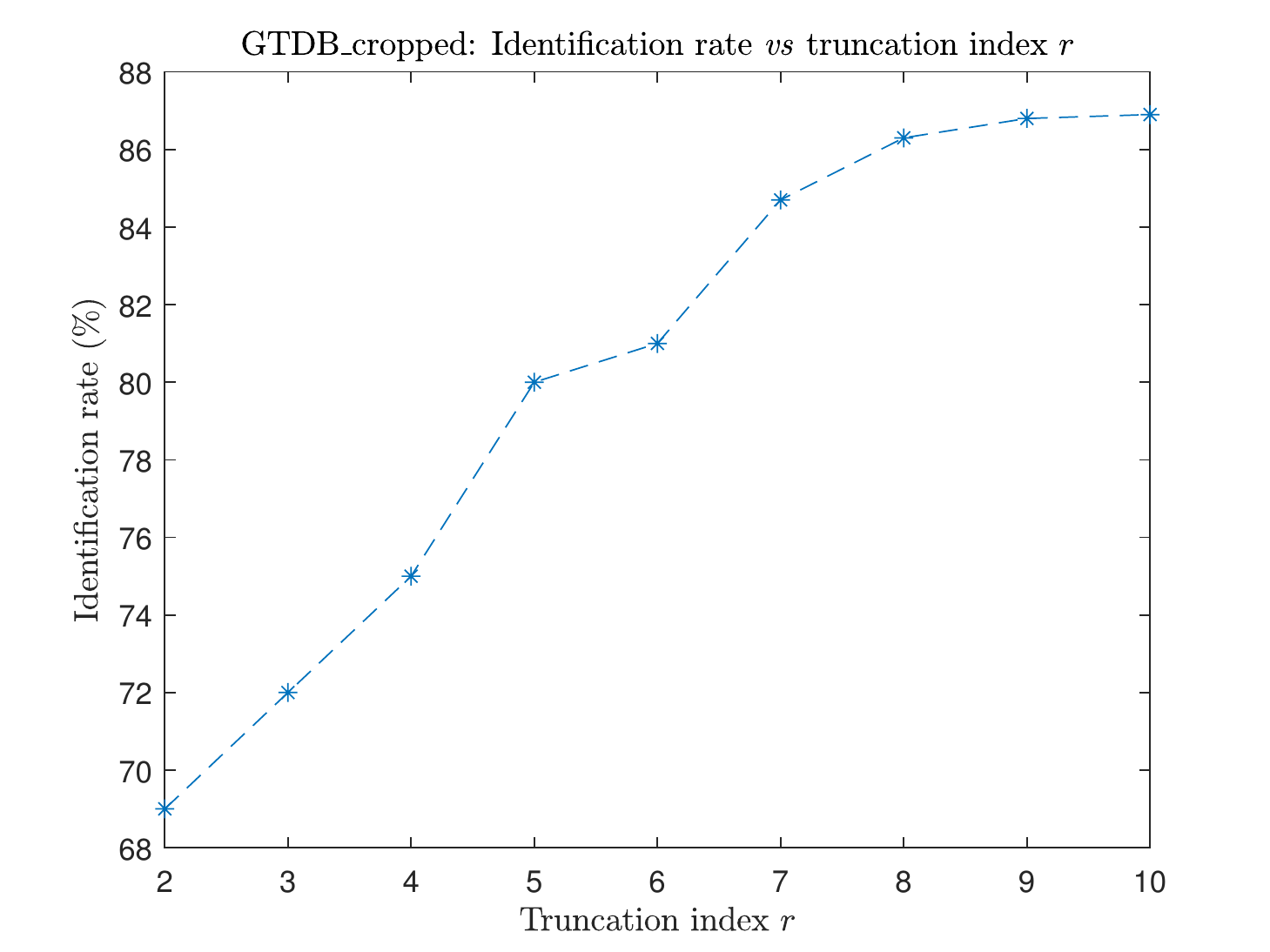}
   \caption{Identification rates for different truncation indices $r$}
   \label{Fig7}
\end{figure}

To compute the rate of recognition, we ran 100 simulations, obtained the number of successes (\textit{ie} a test is successful if the person is correctly identified) and reported the best identification rates, in function of the truncation index $r$ in Fig. \ref{Fig7}.

\noindent The results match the performances observed in the literature \cite{Wang} for this database and it confirms that the use of a Golub Kahan strategy is interesting especially because, in terms of training, the Tube Tensor PCA algorithm required only \textbf{5 seconds} instead of \textbf{25 seconds} when using a c-SVD.  
 
\subsection{Example 3}

In the second example, we used the larger AR face database (cropped version)   \href{ http://cbcsl.ece.ohio-state.edu/protected-dir/AR_warp_zip.zip}{(Face crops)}, \cite{Martinez}, which contains 2600 bitmap pictures of human faces (50 males and 50 females, 26 pictures per person), with different expressions, lightning conditions, facial expressions and face orientation. The bitmap pictures were resized to 100x100 Jpeg images. The same protocol as for Example 1 was followed: we partitioned the set of images in two subsets. Out of 26 pictures, 6 pictures were randomly chosen as test images and the remaining 20 were put into the training folder.

\medskip
\begin{figure}[h!]
   \centering
   \includegraphics[width=9cm]{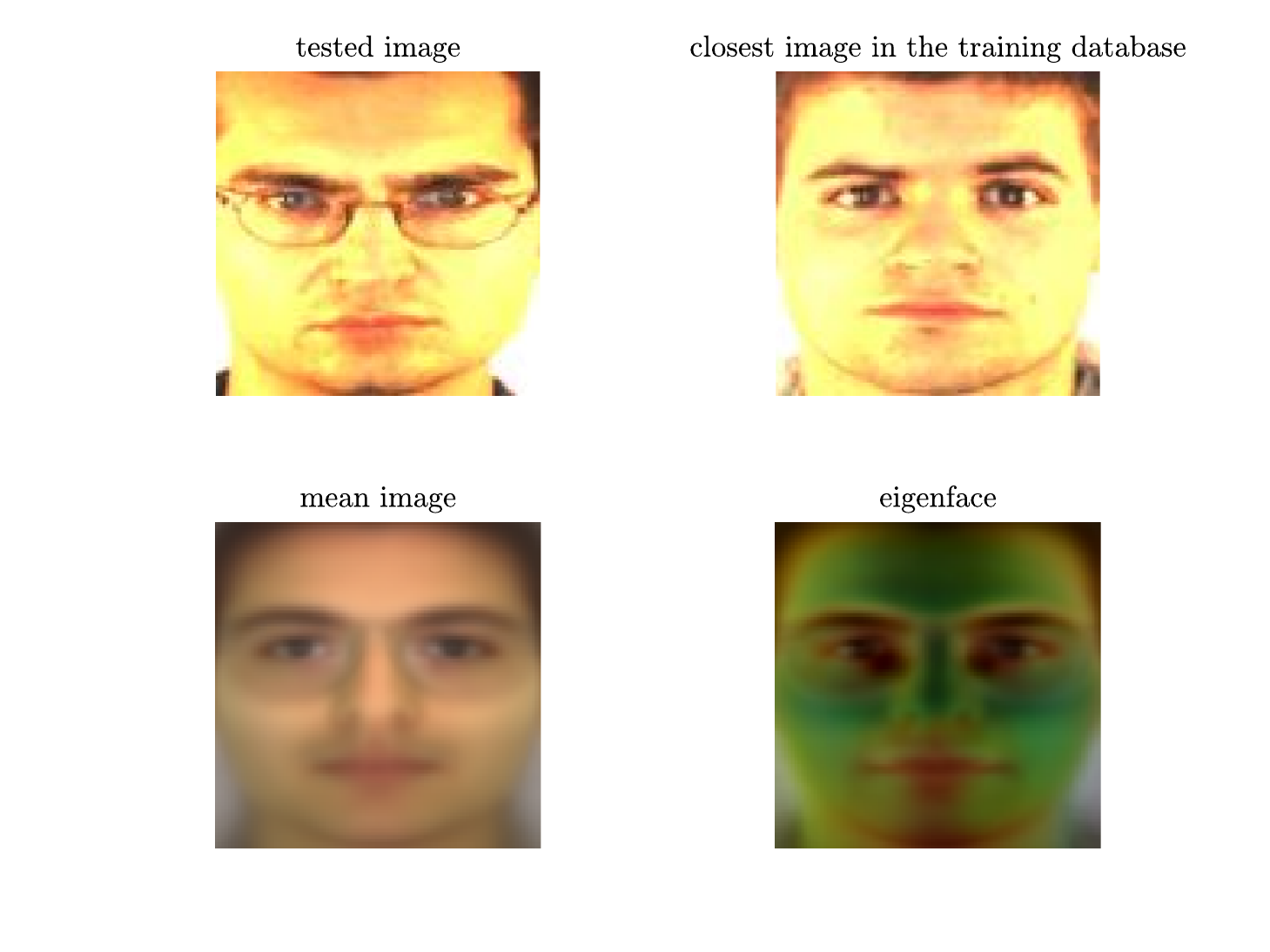}
   \caption{Test image, closest image, mean image and eigenface}
   \label{Fig8}
\end{figure}
\medskip

\begin{figure}[h!]
   \centering
   \includegraphics[width=9cm]{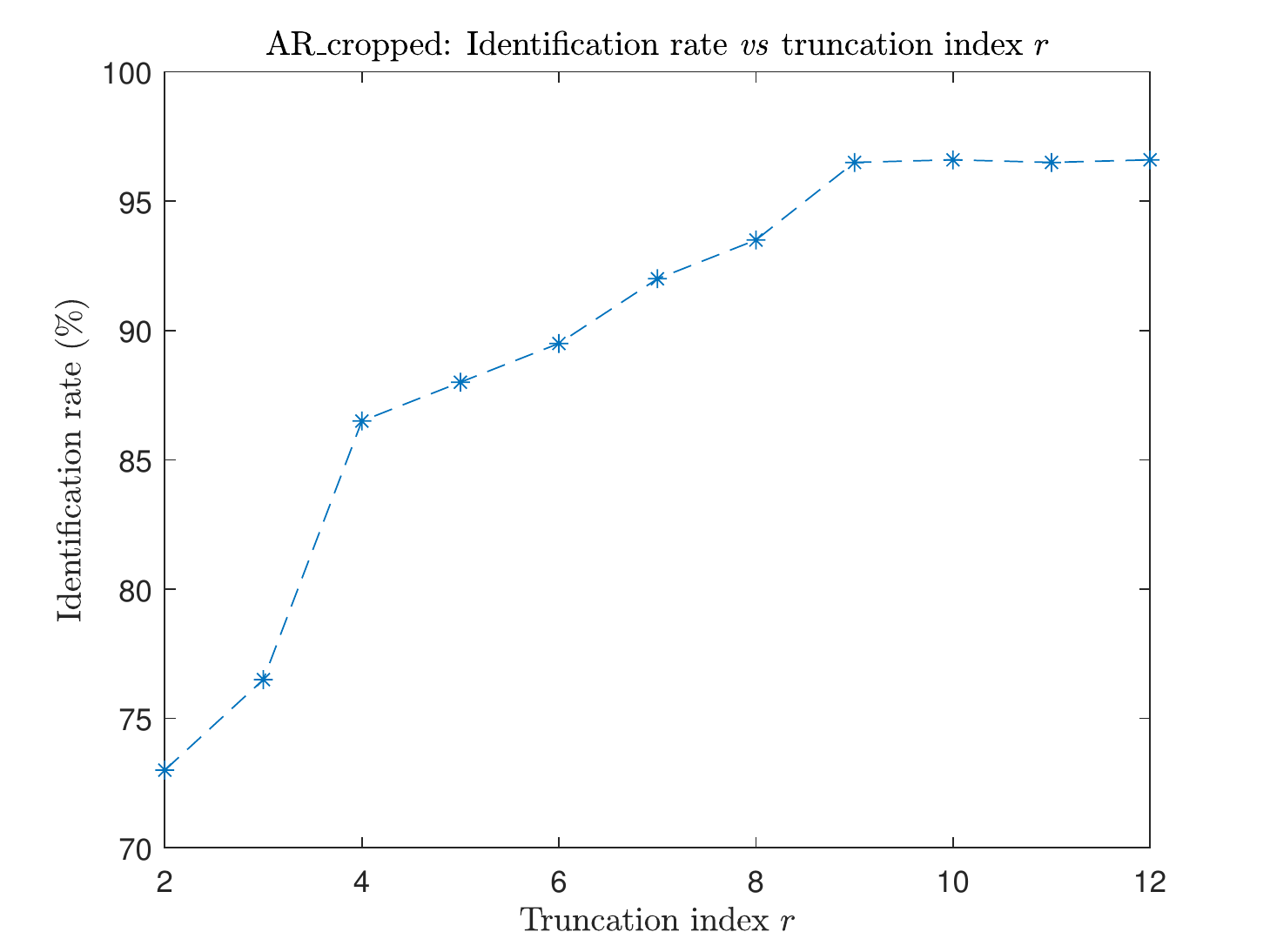}
   \caption{Identification rates for different truncation indices $r$}
   \label{Fig9}
\end{figure}

\noindent 
We applied our approach (TTPCA) to the $10000\times 2000  \times 3$ training tensor $\mathcal{X}$ and plotted the recognition rate as a function of the truncation index in Figure \ref{Fig9}. The training process took \textbf{24 seconds} while it would have taken \textbf{81.5 seconds} if using a c-SVD. \newpage

\noindent For all examples, it is worth noticing that, as expected in face identification problems, only a few of the first largest singular elements suffice to capture the main features of an image. Therefore, the Golub Kahan based strategies such as the TTPCA method are an interesting choice.

\section{Conclusion}

In this manuscript, we focused on two types of Golub Kahan factorizations. We used the recent advances in the field of tensor factorization and showed that this approach is efficient for image identification. The main feature of this approach resides in the ability of the Global Golub Kahan algorithms to approximate the dominant singular elements of a training matrix or  tensor without needing to compute the SVD. This is particularly important as the matrices and tensors involved in this type of application can be very large. Moreover, in the case for which color has to be taken into account, this approach  do not involve a conversion to grayscale, which can be very important for some applications.

\bibliographystyle{plain}

\end{document}